\documentstyle[12pt, fullpage]{amsart}
\begin{document}

\newcommand{\createtitle}[2]{\title{#1}\author{Greg Martin}\address{School of Mathematics\\Institute for Advanced Study\\Olden Lane\\Princeton, NJ 08540\\U.S.A.}\email{gerg@@math.ias.edu}\subjclass{#2}\maketitle}
\newcommand{\make}[1]{\label{#1sec}\noindent\input{#1}}
\newcommand{\Mmake}[1]{\label{#1sec}\noindent}


\newtheorem{theorem}{Theorem}
\newtheorem{lemma}[theorem]{Lemma}
\newtheorem{corollary}{Corollary}[theorem]
\newtheorem{proposition}[theorem]{Proposition}

\newenvironment{pflike}[1]{\noindent{\bf #1}}{\vskip10pt} 
\newenvironment{proof}{\begin{pflike}{Proof:}}{\qed\end{pflike}}

\def\(#1\)#2{
 \def\helpera{\ifcase#2(\or\big(\or\Big(\or\bigg(\or\Bigg(\else(\fi}
 \def\helperb{\ifcase#2)\or\big)\or\Big)\or\bigg)\or\Bigg)\else)\fi}
 \helpera{#1}\helperb}

\newcommand{\3}[1]{\#\{#1\}}
\newcommand{\abs}[1]{\left|#1\right|}
\newcommand{\floor}[1]{\left\lfloor#1\right\rfloor}
\newcommand{\ceil}[1]{\left\lceil#1\right\rceil}

\renewcommand{\mod}[1]{{\ifmmode\text{\rm\ (mod~$#1$)}
 \else\discretionary{}{}{\hbox{ }}\rm(mod~$#1$)\fi}}
\newcommand{\ep}{\varepsilon}
\renewcommand{\implies}{\Rightarrow}
\newcommand{\rmif}{{\rm if\ }}

\newcommand{\half}{{\mathchoice{\textstyle\frac12}{1/2}{1/2}{1/2}}}
\newsymbol\dnd 232D 
\newcommand{\exdiv}{\mathrel{\mid\mid}}

\renewcommand{\lg}[1]{\mathop{\log_{#1}}}
\def\lgs#1^#2{\mathop{\log_{#1}^{#2}}}
\newcommand{\li}{\mathop{\rm li}}

\newcommand{\doublespace}{
  \baselineskip=24pt}
\newcommand{\spaceandahalf}{\parskip6pt\baselineskip=18pt}

\newcommand{\scroll}[1]{\scrollmode #1\errorstopmode}
\newcommand{\comment}[1]{}

\createtitle{Solubility of Systems of Quadratic Forms}{11D72}

\def\v#1{u_F(#1)}
\def\vv#1#2{u_F^{(#1)}(#2)}
\def\qp{{{\bf Q}_p}}
\def\vqp#1{u_\qp(#1)}

It has been known since the last century that a single quadratic form
in at least five variables has a nontrivial zero in any $p$-adic
field, but the analogous question for systems of quadratic forms
remains unanswered. It is plausible that the number of variables
required for solubility of a system of quadratic forms simply is
proportional to the number of forms; however, the best result to date,
from an elementary argument of Leep \cite{Leep:SoQF}, is that the
number of variables needed is at most a quadratic function of the
number of forms. The purpose of this paper is to show how these
elementary arguments can be used, in a certain class of fields
including the $p$-adic fields, to refine the upper bound for the
number of variables needed to guarantee solubility of systems of
quadratic forms. This result partially addresses Problem 6 of Lewis'
survey article \cite{Lew:DPSaU} on Diophantine problems.

By a nontrivial zero of a system of forms $f_1,\ldots,f_t\in
F[x_1,\ldots,x_n]$, we mean a nonzero element {\bf a} of $F^n$ such
that $f_j({\bf a})=0$ simultaneously for $1\le j\le t$. We let $\v
t$ denote the supremum of those positive integers $n$ for which there
exist $t$ quadratic forms over $F$ in $n$ variables with no nontrivial
zero. In other words, assuming $\v t<\infty$, any set of $t$ quadratic
forms in $F[x_1,\ldots,x_n]$, with $n>\v t$, will have a nontrivial
zero (equivalently, a projective zero, since the forms are
homogeneous), while this property does not hold for $n=\v t$. We may
now state our main theorem.

\begin{theorem}Let $F$ be a field, and suppose that for some positive integer
$m$, we have
\begin{equation}
\v m=m\v1.
\label{hypo}
\end{equation}
Then
\begin{equation}
\v t\le\half \big( t(t-m+2)+\tau(m-\tau) \big) \v1,
\label{mainthmeq}
\end{equation}
where $\tau$ is the unique integer satisfying $1\le\tau\le m$ and
$\tau\equiv t\mod m$.
\label{mainthm}
\end{theorem}

We remark that for any $1\le r\le t$, we always have the lower bound
\begin{equation}
\v t\ge \v r+\v{t-r},
\label{easylow}
\end{equation}
for if $f_i(x_1,\ldots,x_{\v r})$ $(1\le i\le r)$ and
$g_j(y_1,\ldots,y_{\v{t-r}})$ $(1\le j\le t-r)$ are systems of
quadratic forms with no nontrivial zeros, then we can combine the two
systems and the two sets of variables to yield a system of $t$
quadratic forms in $\v r+\v{t-r}$ variables with no nontrivial
zeros. In particular, equation (\ref{easylow}) readily implies that
for all $t\ge 1$, we have
\begin{equation}
\v t\ge t\v1\label{bestposs}.
\end{equation}
Thus the hypothesis (\ref{hypo}) of Theorem \ref{mainthm} is a natural one,
representing the best-possible situation for systems of $m$ quadratic
forms.

In fact, if $F$ is a local field (a finite extension either of $\qp$
for some prime $p$, or of $k((T))$ for some finite field $k$), Hasse
\cite{Has:UdDvZetc} has shown that $\v1=4$ (see Lam \cite{Lam:tAToQF}
for an exposition), and Demjanov \cite{Dem:PoQFoetc} has shown that
$\v2=8$ (a simpler proof has been provided by Birch, Lewis, and Murphy
\cite{BirLewMur:SQF}). Thus the following corollary of Theorem
\ref{mainthm} is immediate.

\begin{corollary}Let $F$ be a local field. Then
\begin{equation*}
\v t\le
\begin{cases}
2t^2+2,&t\text{ \rm odd;}\cr
2t^2,&t\text{ \rm even.}
\end{cases}
\end{equation*}
\label{localF}
\end{corollary}

It has also been shown by Birch and Lewis \cite{BirLew:SoTQF}, with a
correction and refinement by Schuur \cite{Sch:oSoTQF}, that whenever
$p\ge11$, we have $\vqp3=12$. Therefore we can again apply
Theorem \ref{mainthm} to obtain the following corollary, which is
superior to Corollary \ref{localF} for these primes.

\begin{corollary}Let $p\ge11$ be prime. Then 
\begin{equation}
u_\qp(t) \le
\begin{cases}
2t^2-2t+4,&t\not\equiv0\mod3;\cr
2t^2-2t,&t\equiv0\mod3.
\end{cases}
\end{equation}
\end{corollary}

The methods employed in this paper are a modest refinement of those of
Leep \cite{Leep:SoQF}, who has shown that $\v t\le\half t(t+1)\v1$ for
arbitrary fields $F$, and also that $\vqp t\le 2t^2+2t-4$ (for $t\ge2$)
for every prime $p$. Because the argument is brief and completely
elementary, we may provide an essentially self-contained proof of
Theorem \ref{mainthm}.

It is a pleasure to thank Trevor Wooley and Hugh Montgomery for their
suggestions on improving this paper and for their guidance in
general. This material is based upon work supported under a National
Science Foundation Graduate Research Fellowship.

\section{Preliminary Lemmas}

Let $\vv dt$ denote the supremum of those positive integers $n$ for
which there exist $t$ quadratic forms over $F$ in $n$ variables whose
set of solutions contain no $(d+1)$-dimensional subspace of $F^n$. In
other words, any set of $t$ quadratic forms in $F[x_1,\ldots,x_n]$,
with $n>\vv dt$, will have a $(d+1)$-dimensional subspace of
simultaneous zeros (or, equivalently, a $d$-dimensional subspace of
projective zeros), while this property does not hold for $n=\vv
dt$. For instance, we have $\vv0t=\v t$.

The following two lemmas can be found in Leep \cite{Leep:SoQF}; we
provide proofs for the sake of completeness.

\begin{lemma}For any field $F$, and for all positive integers $k<t$,
we have
\begin{equation*}
\v t\le \vv{\v k}{t-k}.
\end{equation*}
\label{disection}
\end{lemma}

\begin{proof}
Let $n>\vv{\v k}{t-k}$, and let $f_1,\ldots,f_t$ be quadratic
forms over $F$ in $n$ variables. To establish the lemma, it suffices
to show that these forms have a nontrivial zero in $F^n$. By the
definition of $\vv{\v k}{t-k}$, the system $f_1,\ldots f_{t-k}$ of
$t-k$ quadratic forms has a $(\v k+1)$-dimensional subspace $S$ of
zeros. By parametrizing $S$ with variables $y_1,\ldots,y_{\v k+1}$,
we may consider the restrictions of the forms $f_{t-k+1},\ldots,f_t$
to $S$ as quadratic forms in $\v k+1$ variables. Now by the definition
of $\v k$, these forms have a nontrivial zero in $S$, and so the forms
$f_1,\ldots,f_t$ have a nontrivial zero in $F^n$.
\end{proof}

\begin{lemma}For any field $F$, and for all positive integers $t$ and
$d$, we have
\begin{equation*}
\vv dt \le \vv{d-1}t+t+1.
\end{equation*}
\label{dimshave}
\end{lemma}
\begin{proof}
Let $n>\vv{d-1}t+t+1$, and let $f_1,\ldots,f_t$ be quadratic forms
over $F$ in $n$ variables. To establish the lemma, it suffices to show
that $F^n$ contains a $(d+1)$-dimensional subspace of zeros for these
forms. Since $n>\vv{d-1}t \ge \v t$, we can certainly find a
nontrivial zero for the forms $f_1,\ldots,f_t$, which generates a
1-dimensional subspace $T$ of zeros of these forms. By making a
linear change of variables, we may assume that $T$ is spanned by the
vector $(0,\ldots,0,1)$. For each $1\le j\le t$, we may write
\begin{equation}
f_j(x_1,\ldots,x_n)=x_n^2f_j(0,\ldots,0,1) +
x_nL_j(x_1,\ldots,x_{n-1}) + Q_j(x_1,\ldots,x_{n-1}), \label{factorxn}
\end{equation}
where the $L_j$ and $Q_j$ are linear and quadratic forms,
respectively, in $n-1$ variables (here we are identifying $T^\perp$
with $F^{n-1}$). But we are under the assumption that each
$f_j(0,\ldots,0,1)$ equals 0, and elementary linear algebra allows us
to find a subspace $S$ of $F^{n-1}$ of codimension $t$ on which the
$t$ linear forms $L_1,\ldots,L_t$ all vanish identically. Again we
parametrize $S$ by variables $y_1,\ldots,y_{n-t-1}$ and consider the
restrictions of the forms $Q_1,\ldots,Q_t$ to $S$ as quadratic forms
in $n-t-1>\vv{d-1}t$ variables. By the definition of $\vv{d-1}t$, we
may find a $d$-dimensional subspace $U$ of $S$ consisting of zeros of
the forms $Q_1,\ldots,Q_t$. We now see from (\ref{factorxn}) that
$U\oplus T$ is a $(d+1)$-dimensional subspace of zeros of the original
forms $f_1,\ldots,f_t$.
\end{proof}

\section{Proof of Theorem \ref{mainthm}}

We begin by making some remarks that hold in any field $F$, without
the hypothesis (\ref{hypo}) of Theorem \ref{mainthm}. Using Lemma
\ref{disection} together with several applications of Lemma
\ref{dimshave}, we see that
\begin{equation*}
\v t\le \vv{\v k}{t-k}\le \v{t-k}+(t-k+1)\v k.
\end{equation*}
Therefore, for any positive integer $r$ such that $rk<t$, we have
\begin{equation}
\v t\le \v{t-rk} + \sum_{i=1}^r (t-ik+1)\v k.
\label{laststop}
\end{equation}
Thus we have established a bound for $\v t$ in terms of $\v j$ for
small values of $j$. In fact this is precisely the approach in Leep
\cite{Leep:SoQF}, with the choices $k=1$ and $r=t-1$, so that the
final bound is in terms of $\v1$ alone. One can also choose $r=t-2$
and obtain a bound for $\v t$ in terms of $\v1$ and $\v2$, which will
be better if the value of $\v 2$ is known to be small.

However, for fields $F$ that satisfy the hypothesis (\ref{hypo}) for
some positive integer $m$, it turns out to be more beneficial to take
$k=m$ in the bound (\ref{laststop}). We choose $r$ to make $t-rk$ as
small as possible while still positive: if we let $\tau$ be the
integer satisfying $1\le\tau\le m$ and $\tau\equiv t\mod m$, then
$r=(t-\tau)/m$. With these choices, equation (\ref{laststop}) becomes
\begin{equation}
\v t\le \v\tau + {t-\tau\over2m}(t-m+\tau+2)\v m.
\label{inductwo}
\end{equation}
We claim that $\v m=m\v1$ forces $\v\tau=\tau\v1$ as well, since by
the lower bounds (\ref{easylow}) and (\ref{bestposs}), we have
\begin{equation*}
\begin{split}
\tau\v1 \le \v\tau &\le \v m-\v{m-\tau} \\
&\le m\v1-(m-\tau)\v1=\tau\v1.
\end{split}
\end{equation*}
Substituting these expressions in the bound (\ref{inductwo}) gives us
\begin{equation*}
\v t\le \tau\v1+{t-\tau\over2m}(t-m+\tau+2)m\v1,
\end{equation*}
which is the same as the bound (\ref{mainthmeq}). This establishes the
theorem.

\providecommand{\bysame}{\leavevmode\hbox to3em{\hrulefill}\thinspace}


\begin{thebibliography}{1}

\bibitem{BirLew:SoTQF}
B.~J. Birch and D.~J. Lewis, \emph{Systems of three quadratic forms}, Acta
  Arith. \textbf{10} (1964--1965), 423--442.

\bibitem{BirLewMur:SQF}
B.~J. Birch, D.~J. Lewis, and T.~G. Murphy, \emph{Simultaneous quadratic
  forms}, Amer. J. Math. \textbf{84} (1962), 110--115.

\bibitem{Dem:PoQFoetc}
V.~B. Demjanov, \emph{Pairs of quadratic forms over a complete field with
  discrete norm with a finite field of residue classes}, Izv. Akad. Nauk SSSR,
  Ser. Mat. \textbf{20} (1956), 307--324.

\bibitem{Has:UdDvZetc}
H.~Hasse, \emph{{\"Uber} die {D}arstellbarkeit von {Z}ahlen durch quadratische
  {F}ormen im {K\"orper} der rationalen {Z}ahlen}, J. Reine Angew. Math.
  \textbf{152} (1923), 129--148.

\bibitem{Lam:tAToQF}
T.~Y. Lam, \emph{The algebraic theory of quadratic forms}, Benjamin/Cummings
  Publishing Co., Inc., Reading, Mass., 1973.

\bibitem{Leep:SoQF}
D.~B. Leep, \emph{Systems of quadratic forms}, J. Reine Angew. Math.
  \textbf{350} (1984), 109--116.

\bibitem{Lew:DPSaU}
D.~J. Lewis, \emph{Diophantine problems: Solved and unsolved}, Number Theory
  and Applications (R.~A. Mollin, ed.), Kluwer Academic Publishers, Dordrecht,
  1989, pp.~103--121.

\bibitem{Sch:oSoTQF}
S.~E. Schuur, \emph{On systems of three quadratic forms}, Acta Arith.
  \textbf{36} (1980), 315--322.

\end{thebibliography}
\end{document}